\documentclass{elsart}
\usepackage{graphicx,amsmath, amsfonts, amssymb, array}
\usepackage[english]{babel}
  \usepackage{paralist}
  \usepackage{graphics} %% add this and next lines if pictures should be in esp format
  \usepackage{epsfig} %For pictures: screened artwork should be set up with an 85 or 100 line screen
\newtheorem{theorem}{Theorem}[section]
\newtheorem{corollary}{Corollary}
\newtheorem{lemma}[theorem]{Lemma}

\theoremstyle{definition}
\newtheorem{definition}[theorem]{Definition}
\newtheorem{remark}{Remark}
\newtheorem{example}[corollary]{Example}
\newcommand{\tr}{{\rm tr}\,}

\newcommand{\Tr}{\mathrm{Tr}}

\newcommand{\R}{\mathbb{R}}
\newcommand{\C}{\mathbb{C}}
\newcommand{\F}{\mathbb{F}}
\newcommand{\Z}{\mathbb{Z}}
\newcommand{\Q}{\mathbb{Q}}

\DeclareMathOperator{\GF}{GF}

\DeclareMathOperator{\soc}{Soc}
\DeclareMathOperator{\rad}{Rad}

\begin{document}
\begin{frontmatter}

\title
{On the Weight Distribution of Codes over Finite Rings}
\author{Eimear Byrne}
\address{School of Mathematical Sciences, University College Dublin}
 \ead{ebyrne@ucd.ie}
 \thanks{Research supported by the Science Foundation of Ireland, 08/RFP/MTH1181}
%\subjclass{Primary: 11T71; Secondary: 14G50}


\begin{abstract}
   Let $R>S$ be finite Frobenius rings for which there exists a trace map
   $T:{_S}R \longrightarrow {_S}R $. Let 
   $C_{f,S}:=\{x \mapsto T(\alpha x + \beta f(x)) : \alpha, \beta \in R \}$.
   $C_{f,S}$ is an $S$-linear subring-subcode of a left linear code over $R$.
   We consider functions $f$ for which the homogeneous weight distribution
   of $C_{f,S}$ can be computed. In particular, we give constructions 
   of codes over integer modular rings and commutative local Frobenius that have small spectra.
\end{abstract}

\end{frontmatter}

\section{Introduction}

        The homogeneous weight, introduced for integer residue rings in \cite{ch} and 
        extended for 
        arbitrary finite rings in \cite{bruce}, has been studied extensively in the 
        context of ring-linear coding. 
        It can be viewed as a generalization of the Hamming weight; in fact it 
        coincides with the Hamming weight when the underlying ring is a finite field 
        and is the Lee weight when defined over $\mathbb{Z}_4$. 
        Many of the classical results for codes over finite 
        fields for the Hamming weight have corresponding homogeneous weight versions 
        for codes over finite rings. The 
        MacWilliams equivalence theorem holds for codes over finite Frobenius rings and 
        quasi-Frobenius modules with respect to the homogeneous weight 
        \cite{bruce,gnw}. Analogues of several classical bounds have been 
        found  for this weight function \cite{bgos,bgks,plotkin}. Combinatorial objects 
        such as strongly regular graphs can be constructed from codes over finite 
        Frobenius rings with exactly two nonzero homogeneous weights \cite{bgh,bs,h07}.  
        Although homogeneous weights exist on any finite ring, if the ring in question 
        is Frobenius, these weight functions can be expressed in terms of a character 
        sum \cite{h}, a property we shall use here.        
        
        Let $\Tr:\GF(q^r) \longrightarrow GF(q)$ be the usual trace map from $GF(q^r)$ 
        onto $GF(q)$. Let $f:\GF(q^r) \longrightarrow GF(q^r)$ be an arbitrary map.
        A construction of a $GF(q)$-linear subspace-subcode is given by:
        $$C_f:=\{c^f_{\alpha,\beta}:\GF(q^r) \longrightarrow GF(q):x \mapsto \Tr(\alpha 
        x+\beta f(x)), \alpha,\beta \in GF(q^r)\}.$$
        This has arisen in the literature on perfect and almost perfect nonlinear 
        functions (c.f. \cite{bbmmg,CCZ,ycd}) and on cyclic codes (in which case $f$ is 
        a power map). 
        In the case of an APN function on a field of 
        even characteristic it is sometimes possible to determine the distinct weights 
        or the weight distribution of the resulting code, or equivalently the Walsh 
        spectrum of $f$. In \cite{cdy}, perfect nonlinear maps were used to construct 
        (in some cases optimal) codes for use in secret sharing schemes.
        
        Here we consider the same code construction for codes over finite rings, 
        specifically, some integer modular rings, Galois rings and local commutative 
        Frobenius rings. 
        
\section{Preliminaries}
 
         We recall some properties of finite rings that meet our purposes. 
         Further details can be read in \cite{h,lam,macd,r}.
         For a finite ring $R$, we denote by $\hat{R}:=$ Hom$_{\Z}(R,\C^{\times})$,
         the group of additive characters of $R$. $\hat{R}$
         is an $R$-$R$ bi-module according to the relations
         $${^r\chi}(x) = \chi(rx),\:\:\: \chi^r(x) = \chi(xr),$$ 
         for all $x,r \in R, \chi \in \hat{R}$. A character $\chi$ is called left (resp. right) generating
         if given any $\phi \in \hat{R}$ there is some $r \in R$ satisfying 
         $\phi = {^r}{\chi}$ (resp. $\phi = {\chi}^r$). 
         The next result gives a characterization of finite Frobenius rings.    
         \begin{theorem}
         Let $R$ be a finite ring. The following are equivalent.
         \begin{enumerate}
         \item
         $R$ is a Frobenius ring 
         \item
         $\soc {_R}R$ is left principal,
         \item
         ${_R}(R/Rad\; R) \simeq  \soc{_R}R$, 
         \item
         $ {_R}R \simeq {_R}\hat{R} $ 
         \end{enumerate}
         \end{theorem}

         Then $~_R\hat{R}= ~_R\langle \chi \rangle$ for some (left) generating character $\chi$. It can be shown that any left generating character is also a right generating character.
 
 Integer residue rings, finite chain rings, semi-simple rings, principal ideal rings, direct products of Frobenius rings, matrix rings over Frobenius rings, group rings over Frobenius rings are all examples of Frobenius rings. The results of this paper are restricted to the case where the code's alphabet is an integer modular ring or a local commutative Frobenius ring.

 Let $R$ be a finite commutative local Frobenius ring with unique maximal ideal $M$ and residue field $K=R / M$  
 of order $q$ for some prime power $q$. Then $M = \rad R$ and $\soc R = R / M$ is simple. 
 Moreover, $\soc R$ is the annihilator  of $\rad R$, which we write as $\soc R = M^\perp$.
 $R^\times$ contains a unique cyclic subgroup $G$ of order $|K^\times|$. We call ${\mathcal T}:=G \cup \{0\}$ the 
 {\em Teichmuller set} of $R$.
 Later, we will use that fact that each element $a \in R$ can be expressed uniquely as
 $a = a_t + a_m $ for some unique $a_t \in {\mathcal T}$, $a_m \in M$. 
 We define the map $\nu: R \longrightarrow {\mathcal T}: a \mapsto a_t$.
 
 If the local commutative ring $R$ is the Galois ring 
 $GR(p^{n},r)$, of order $p^{nr}$ and characteristic $p^n$, then
 each element $a \in R$ can be expressed uniquely in the form
 $a=a_0 + pa_1+\cdots +p^{n-1}a_{n-1}$ for some unique $a_i \in {\mathcal T}$.

  For an arbitrary finite ring, the homogeneous weight is defined as follows \cite{ch,bruce}.

  \begin{definition}
  Let $R$ be a finite ring. A weight $w:R\longrightarrow \Q$ is
   \emph{(left) homogeneous}, if $w(0)=0$ and
    \begin{enumerate}
  \item If $Rx=Ry$ then $w(x)= w(y)$ for all $x,y\in R$.
  \item There exists a real number $\gamma$ such that
    \begin{equation*}
  \sum_{y\in Rx}w(y) \; =\; \gamma \, |Rx|\qquad\text{for all $x\in
    R\setminus \{0\}$}.
\end{equation*}
\end{enumerate}
\end{definition}
 
 \begin{example}
  On every finite field $\F_q$ the Hamming weight
    is a homogeneous weight of average value $\gamma=\frac{q-1}{q}$.
 \end{example}

\begin{example}
On the ring $\Z_{pq}, p,q$ prime, a homogeneous weight with average value $\gamma = 1$ is given by
 \begin{equation*}
      w: R\longrightarrow {\mathbb R},\quad x \mapsto
      \left\{\begin{array}{cl}
                     0& \mbox{if } x=0,\\
        \frac{p}{p-1} & \mbox{if } x \in p\Z_{pq} ,\\
        \frac{q}{q-1} & \mbox{if } x \in q\Z_{pq},\\
        \frac{pq-p-q}{pq-p-q+1} & {\mbox{ otherwise.}}       
      \end{array}\right.
    \end{equation*}
\end{example}

\begin{example}
On a local Frobenius ring $R$ with $q$-element residue field the
    weight
    \begin{equation*}
      w: R\longrightarrow {\mathbb R},\quad x \mapsto
      \left\{\begin{array}{cl}
        0            & \mbox{if } x=0,\\
        \frac{q}{q-1}& \mbox{if } x\in \soc(R),\; x\neq 0,\\
        1            & \mbox{if }  \mbox{otherwise},
      \end{array}\right.
    \end{equation*}
is a homogeneous weight of average value $\gamma=1$.
\end{example}

%\end{document}
A description of the homogeneous weight in terms of sums of generating characters is given by the following \cite{h}.
\begin{theorem} 
Let $R$ be a finite Frobenius
  ring with generating character $\chi$. 
  Then the homogeneous weights on $R$ are precisely the functions
  \begin{equation*}
w: R\longrightarrow {\mathbb R}, \quad x \mapsto
  \gamma\Big[1-\frac{1}{|R^{\times}|}\sum_{u\in R^{\times}}
  \chi(xu)\Big]
\end{equation*}
where $\gamma$ is a real number.
\end{theorem}
Unless otherwise stated, we will set $\gamma=1$ (the {\em normalized} homogeneous weight).

\section{Characters and Trace Maps} 

\begin{definition}
Let $R>S$ be Frobenius rings.  
An $S$-module epimorphism $T:{_S}R\longrightarrow {_S}{S}$ 
whose kernel contains no non-trivial left ideal of $R$ is called a trace map from $R$ onto $S$.
\end{definition}

\begin{example}
Recall that a finite Frobenius ring $R$ has a generating character $\chi$. Let $R$ have 
characteristic $m$. Then we can implicitly define a trace map $T$ from $R$ onto its characteristic subring $\Z_m$ by
$\chi(x) = \omega^{T(x)}$ for all $x$, where $\omega$ is a primitive complex $m$th root of unity.
This is the absolute trace map on $R$. 
\end{example}

Given a trace map $T:{_S}R\longrightarrow {_S}{S}$, a generating character $\Phi \in \hat{S}$ determines a generating character
$\chi \in \hat{R}$ by:
$$\chi(x) = \Phi(T(x)) \; \forall \; x \in R. $$

\begin{example}
Let $S$ be a finite Frobenius ring and let $R=M_n(S)$. Then $R$ is Frobenius and the usual trace map
$$\tr:R \longrightarrow S : (a_{ij}) \mapsto \sum_{i=1}^n a_{ii}$$
is an epimorphism onto $S$ whose kernel contains no non-trivial left ideal of $R$.
A generating character $\Phi \in \hat{S}$ induces a generating character $\chi = \tr \circ \Phi  \in \hat{R}$.
\end{example}

\begin{example}\label{exgaltr}
Let $R=GR(p^n,sk), S:=GR(p^n,s)$ be Galois rings of characteristic $p^n$ and orders $p^{nsk},p^{ns}$, respectively.
As in the case of a finite field, $R$ has a cyclic automorphism group of order $sk$ \cite{r}.
Each element $a\in R$ has a canonical representation in the form $a = \sum_{i=0}^n p^i a_i $ for some unique $a_i$ in the Teichmuller set of $R$. 
With respect to this expression,
$$\sigma: R \longrightarrow R :\sum_{i=0}^n p^i a_i \mapsto  \sum_{i=0}^n p^i a_i^{p}$$
generates $Aut(R)$, and $\tau:=\sigma^s$ generates $Aut(R:S)$, the group of $S$-automorphisms of $R$.
The map
$$T_{R/S}: R \longrightarrow S : a \mapsto a + \tau(a) + \cdots + \tau^{(k-1)}(a)$$
is a trace map from $R$ onto $S$.  Observe that as in the field case, $T_{R/S}(\tau(a)) = T_{R/S}(a)$ for any $a \in R$.
\end{example}

\begin{example}
    Let $R=\Z_4[x]/\langle x^2+2 \rangle$. Then every element $r \in R$ can be expressed uniquely in the form $r = r_0+\theta r_1$, where $\theta^2 = 2$ $r_i \in \Z_4$. 
A $\Z_4$-epimorphism from $R$ onto $\Z_4$ must have the form
$T_\lambda: R \longrightarrow \Z_4 : r_0+\theta r_1 \mapsto \lambda_0 r_0 + \lambda_1 r_1$  for some $\lambda_i \in \Z_4$.
This gives a trace map if and only if $\lambda_1 \in \Z_4^\times$, so there are exactly $8$ distinct trace maps from $R$ onto $\Z_4$. Each such map determines a generating character $\chi$ defined by
$\chi^\lambda (x)= \omega^{T_\lambda(x)} $ for all $x \in R.$
Note that $Aut(R)$ is cyclic of order 2 generated by $\sigma: r_0+\theta r_1 \mapsto   r_0 - \theta r_1.$
It is easy to check that $T_\lambda \circ \sigma  = T_{\sigma(\lambda)} \neq T_\lambda$ for 
$\lambda_1 \in \Z_4^\times.$
\end{example}

\section{Subring Subcodes}

For the remainder, let $R>S$ be finite Frobenius rings and assume there is a trace map
$T:{_S}R \longrightarrow {_S}S$. Let $\Phi$ be a generating character of $S$ and let 
$\chi := \Phi \circ T$.
For any map $f:R \longrightarrow R$, we define the left $S$-linear subring subcode
$$C_{f,S} := \{c^f_{\alpha,\beta}: R \longrightarrow S :x \mapsto T(\alpha x + \beta f(x)): \alpha,\beta \in R  \}.$$ 
For the case $R=S$, we write $C_f:=C_{f,R}$. We extend the homogeneous weight on $S$ to the module of functions from $R$ to $S$ by $w(g):= \sum_{x \in R} w(g(x))$ for an arbitrary map $g:R \longrightarrow S$.
We thus compute the weight of each codeword as: 
\begin{eqnarray*}
w(c^{f}_{\alpha,\beta}) &=& \sum_{x \in R} w(c^f_{\alpha,\beta}(x)) \\
 & = & |R| - \frac{1}{|S^\times|} \sum_{u \in S^\times} \sum_{x \in R} \Phi^u(T(\alpha x + \beta f(x)))\\
 & = &|R| - \frac{1}{|S^\times|} \sum_{u \in S^\times} \sum_{x \in R} \chi^u(\alpha x + \beta f(x)).
\end{eqnarray*}

\begin{definition} 
 Let $f:R \longrightarrow R$. For each $\alpha, \beta \in R$, we define the transform 
 $$W^{f,S}(\alpha,\beta): =\frac{1}{|S^\times|} \sum_{u \in S^\times} \sum_{x \in R} \chi^u(\alpha x + \beta f(x)) 
 = |R| - w(c^f_{\alpha,\beta}).$$
 The spectrum of $f$ is defined to be the set
 $$\Lambda_{f,S} := \{ W^{f,S}(\alpha,\beta): \alpha, \beta \in R \}.$$
\end{definition}

In keeping with the above, we write $W^{f}:=W^{f,R}$ and $\Lambda_f:=\Lambda_{f,R}$.

Observe that if $|\Lambda_{f,S}| = k+1$ then $C_{f,S}$ has exactly $k$ non-zero weights.
Clearly,
$$W^{f,S}(\alpha,0) = \left\{ \begin{array}{ll} \displaystyle{\frac{1}{|S^\times|} \sum_{u \in S^\times} \sum_{x \in R} \chi^u(\alpha x) = |R|} & \mbox{if } \alpha = 0\\
 0 & \mbox{if } \alpha \neq 0.
 \end{array} \right.$$  
In particular, if $f$ is 
$R$-linear then $C_{f,S}$ is a one-weight code with constant non-zero weight $|R|$.

For a code $C$ with distinct homogeneous weights $W \subset \Q$ we define the (homogeneous) weight enumerator of $C$ by $\displaystyle{\sum_{w \in W} A_w X^w}$, where $A_w =|\{c \in C : w(c) = w\}|$\footnote{Note that as defined, this weight enumerator is not necessarily a polynomial, although this could be achieved by choosing appropriate $\gamma$, for example, $\gamma = |R^\times|$.}.

\section{Families of Codes with Few Weights}

We now present some families of codes whose weight distributions can be computed.

\subsection{Galois Rings}

The class of functions discussed in the next theorem is a variant of the family of Frank sequences (c.f. \cite{frank}) .
Such sequences and their generalizations form a family of perfect sequences,
having the property that their Fourier transforms have constant absolute value. 

 \begin{theorem}\label{thgalfrank}
    $R=GR(p^2,r), S=GR(p^2,s), p$ prime, for some positive integers $r,s,k$ satisfying $sk=r$ and $k>1$. 
    Let ${\mathcal T}$ denote the Teichmuller set of $R$, let $q=p^s$ and let $\pi$ be a permutation of ${\mathcal T}$ 
    that fixes $0$. Write $x=x_0+px_1$, $x_i \in {\mathcal T}$  for each $x \in R$. Let 
    $$f:R \longrightarrow R : x \mapsto p\pi(x_0)x_1. $$
    Then $$\Lambda_{f,S} = \{\displaystyle{q^{2k}, q^{k}, -\frac{q^k}{q-1},0} \}.$$
 \end{theorem}
 
 {\bf Proof:}
 Let $\alpha,\beta \in R$. Then $\alpha x + \beta f(x) = \alpha x_0 + p (\alpha+\beta \pi(x_0)) x_1,$
 and so
 \begin{eqnarray*}
    W^{f,S}(\alpha,\beta) & = &\frac{1}{|S^\times|} \sum_{u \in S^\times} \sum_{x \in R} \chi^u(\alpha x + \beta f(x)) \\
    & = &  \frac{1}{|S^\times|} \sum_{u \in S^\times} \sum_{x \in R} \chi^u(\alpha x_0 + p (\alpha+\beta \pi(x_0)) x_1) \\
    & = & \frac{1}{|S^\times|} \sum_{u \in S^\times} \sum_{x_0 \in {\mathcal T}}  \chi^u(\alpha x_0)  \sum_{x_1 \in {\mathcal T}} 
    \chi^u(p(\alpha + \beta \pi(x_0))x_1) \\
    & = &  \frac{1}{|S^\times|} \sum_{u \in S^\times}  \sum_{t \in {\mathcal T}}  \chi^u(\alpha t) \sum_{z \in pR}
      \chi^u((\alpha + \beta \pi(t))z) \\
    & =  & \frac{|pR|}{|S^\times|} \sum_{u \in S^\times}  \sum_{t \in {\mathcal V}}  \chi^u(\alpha t), 
 \end{eqnarray*}
 where ${\mathcal V} = \{ t \in {\mathcal T} : \alpha + \beta \pi(t) \in pR \}.$ 
 
 For the case $\beta \in pR$, we have $W^{f,S}(\alpha,\beta) = W^{f,S}(\alpha,0)$, which takes the value zero if $\alpha \neq 0$ and $|R|$ otherwise. 
 %if $\alpha \in R^\times$ then ${\mathcal V}$ is empty and $W^{f,S}(\alpha,\beta) = 0.$
 %If $\alpha \in pR$ then ${\mathcal V} = {\mathcal T}$ and 
 %$$W^{f,S}(\alpha,\beta) = \left\{ \begin{array}{ll} \displaystyle{\frac{|pR|}{|S^\times|} \sum_{u \in S^\times}  \sum_{z \in pR}  \chi^u(x)} = 0 & \mbox{if }\alpha \neq 0 \\
  % \displaystyle{\frac{|pR|}{|S^\times|}|S^\times||pR|=|pR|^2=|R| = q^{2k}} & \mbox{if } \alpha =0   \end{array} \right.$$
 For arbitrary $\alpha,\beta \in R$ we have
  \begin{eqnarray*}
    W^{f,S}(\alpha,\beta) & = &  \frac{|pR|}{|S^\times|}  \sum_{t \in {\mathcal V}} 
    \left(\sum_{u \in S}  \chi(\alpha tu) - \sum_{u \in pS}  \chi({\alpha t}u) \right) \\
    & = &  \frac{|pR|}{|S^\times|}  \sum_{t \in {\mathcal V}} 
    \left(\sum_{u \in S}  \Phi(T_{R/S}(\alpha t)u) - \sum_{u \in pS}  \Phi(T_{R/S}(\alpha t)u) \right) \\
    & = &  \frac{|pR|}{|S^\times|}(|S||{\mathcal V} \cap {U_0}| - |pS||{\mathcal V} \cap {U_1}|),
    \end{eqnarray*}
    where $U_0 =\{ t \in {\mathcal T} : T_{R/S}(\alpha t) = 0 \}$ and $U_1 = \{ t \in {\mathcal T} : T_{R/S}(\alpha t) \in pR \}.$
    
    For the case $\beta \in R^\times$, ${\mathcal V} = \{\pi^{-1}(\nu(-\alpha \beta^{-1})) \}$.
    If $\alpha \in pR$ then $W \cap U_0 = W \cap U_1 = W = \{0 \}$ and so $W^{f,S}(\alpha,\beta) = |pR|$. If $\alpha \in R^{\times}$, then 
    $$W^{f,S}(\alpha,\beta) = \left\{ \begin{array}{ll} 
                                   \displaystyle{\frac{|pR|}{|S^\times|}(|S| - |pS|) = |pR| = q^k} & \mbox{if } T_{R/S}(\alpha \pi^{-1}(\nu(-\alpha\beta^{-1} ) ) )=0,\\
                                   \displaystyle{\frac{|pR|}{|S^\times|}( - |pS|) = -\frac{q^k}{q-1}}& \mbox{if } T_{R/S}(\alpha \pi^{-1}(\nu(-\alpha\beta^{-1}))) \in pS \backslash \{0\},\\
                                   0 &\mbox{otherwise}.                               
                                  \end{array} \right.$$
   The result follows.
     \qed

 \begin{corollary}
    Let $C_{f,S}$ be defined as in Theorem \ref{thgalfrank}. Then
    $C_{f,S}$ has size $q^{3k}$, and weight enumerator
    \begin{eqnarray*}
    1+(q^k-1)\left((q^{2k-2}-q^{k-1}+q^{k})X^{q^{2k}-q^k}+(q^{2k}-q^{2k-1}+q^k+1)X^{q^{2k}}\right. \\
             \left.+(q^{2k-1}-q^k-q^{2k-2}+q^{k-1})X^{q^{2k}+ \frac{q^k}{q-1}}\right).
    \end{eqnarray*}
  \end{corollary}
  
  {\bf Proof:} Since $\ker T_{R/S}$ contains no non-trivial ideal, $T_{R/S}(\alpha x + p\beta \pi(x_0)x_1)=0$ for all $x \in R$ if and only if $\alpha = 0$ and $\beta \in pR$. It follows that $|C_f| = q^{3k}$. Moreover, given any $\alpha,\alpha',\beta,\beta' \in R$, $c^f_{\alpha,\beta} = c^f_{\alpha',\beta'} $ if and only if $\alpha=\alpha'$ and $\beta - \beta' \in pR$.
  For $\beta \in R^\times$, the map $\alpha \mapsto \alpha\pi^{-1}(\nu(-\alpha \beta^{-1}))$ is a permutation on $R^\times$. Then 
  $|\{ \alpha \in R^{\times}:  T_{R/S}(\alpha \nu(\pi(-\alpha \beta^{-1}))) = \theta \}|$ = $|\{ \alpha \in R^{\times}:  
  T_{R/S}(\alpha) = \theta \}|$ for any $\theta \in R$.
  %Now $T_{R/S}(\alpha) \in pS$ if and only if $\alpha \in \ker T_{R/S} + pR$, in which case 
  %$\eta(\alpha) \in \ker \Tr_{GF(q^k)/GF(q)}$. 
  %Now, $|R^\times \cap \ker T_{R/S} + pR| = q^{k-1}(q^k-1)$.
  Now $|\ker T_{R/S} \cap R^{\times}| =|\ker T_{R/S}| - |\ker T_{R/S} \cap pR| = q^{2(k-1)} - q^{k-1}$ and 
  $|(\ker T_{R/S} + pR) \cap R^{\times}| =|\ker T_{R/S}+pR| - |(\ker T_{R/S}+pR) \cap pR| = q^{2k-1} - q^{k}$.   We summarize these observations along with results of Theorem \ref{thgalfrank} in the table shown below, from which the statement of the theorem follows.
  \begin{center}  
  \begin{tabular}{|c|l|l|}
  \hline
  $w(c^f_{\alpha,\beta})$ & constraints on $\alpha,\beta$ & $|\{ c^f_{\alpha,\beta}\}|$\\
  \hline
   0 & $\beta = 0$, $\alpha = 0$. & 1\\
   \hline
   $q^{2k}-q^k$ & $\beta \in {\mathcal T} \backslash \{0\}, \alpha \in pR$; & \\
                & $\beta \in {\mathcal T} \backslash \{0\}, \alpha \in R^{\times}$ and & $(q^{k}-1)q^k+(q^{k}-1)(q^{2(k-1)}-q^{k-1})$\\
                & $T_{R/S}(\alpha \nu(-\alpha \beta^{-1})) = 0.$ 
                &\\  
   \hline              
   $q^{2k}$     & $\beta = 0, \alpha \neq 0$; & \\
                & $\beta \in {\mathcal T} \backslash \{0\}, \alpha \in R^{\times}$ and & $q^{2k}-1+(q^{k}-1)(q^{2k}-q^{2k-1})$\\
                & $T_{R/S}(\alpha \nu(-\alpha \beta^{-1})) \notin pS$. 
                &\\
   \hline             
   $\displaystyle{q^{2k}+\frac{q^k}{q-1}}$  &   $\beta \in {\mathcal T} \backslash \{0\}, \alpha \in R^{\times}$ and  & $(q^k-1)q^{k-1}(q^{k-1}-1) (q-1)$ \\ 
    & $T_{R/S}(\alpha \nu(-\alpha \beta^{-1})) \in pS \backslash \{0\}$. & \\        
  \hline  
  \end{tabular}
  \end{center}
  ~\qed
  
  For the case $R=S$, the code corresponding to $f:x \mapsto p\pi(x_0)x_1$ is a two-weight code. 
  
  \begin{corollary}
      Let $f$ be defined as in Theorem \ref{thgalfrank}.
    Then $$\Lambda_{f} = \{p^{2r},p^{r},0 \},$$
    and $C_{f}$ has weight enumerator
    $$ 1+(p^{2r}-p^r)X^{p^{2r}-p^r} + (p^{3r}-p^{2r}+p^r-1) X^{p^{2r}}. $$
  \end{corollary}
  
  {\bf Proof:} Let $\alpha,\beta \in R$. 
  Almost exactly as in the proof of Theorem \ref{thgalfrank} we have
  $$W^{f}(\alpha,\beta) =  \frac{|pR|}{|R^\times|}(|R||{\mathcal V} \cap {U_0}| - |pR||{\mathcal V} \cap {U_1}|),$$
    where $U_0 =\{ t \in {\mathcal T} : \alpha t = 0 \}$ and $U_1 = \{ t \in {\mathcal T} : \alpha t \in pR \}.$
    If $\alpha \in R^\times$ then $W \cap U_0=W \cap U_1 = \emptyset$.
  \qed

\subsection{Integer Modular Rings}

%For the $\Z_p$ case, the codes determine orthogonal arrays.

Cyclic codes on finite fields have been well studied. It is surely well known and not hard to show the following.

\begin{theorem}\label{thzp}
  Let $R=\Z_p$, $p$ prime, let $d \in \{2,...,p-1\}$ then $C_{x^d}$ is a two-weight code with Hamming weight 
  enumerator
  %$$1 + \frac{(p-1)^2}{\ell}X^{p\frac{p-\ell-1}{p-1}} + (p^2-1 - \frac{(p-1)^2}{\ell})X^{p},$$ 
  $$1 + \frac{(p-1)^2}{\ell}X^{p-\ell-1} + \left(p^2-1 - \frac{(p-1)^2}{\ell}\right)X^{p-1},$$ 
  where $(d-1,p-1)=\ell$.
\end{theorem} 

%{\bf Proof:}
%It is clear that $C_{x^d}$ has $p^2$ codewords. Let $\alpha, \beta \in \Z_p$. %We compute
%%\begin{eqnarray*}
%%W^f(\alpha,\beta) & = & \frac{1}{p-1}\sum_{x \in \Z_p}\left(\sum_{u \in \Z_p} 
%%\chi(u(\alpha x+\beta x^d))-1 \right)
%%                   =  \frac{p}{p-1}\left( |{\mathcal V}_{\alpha,\beta}|  - 1 \right), 
%%\end{eqnarray*}
%The Hamming weight of $c^f_{\alpha,\beta}$ is $p- |{\mathcal V}_{\alpha,\beta}| $
%where ${\mathcal V}_{\alpha,\beta} = \{x \in \Z_p : \alpha x + \beta x^d = 0 \}$. Clearly, ${\mathcal %V}_{\alpha,0}=\{0\}$ if $\alpha \neq 0$ and is $\Z_p$ otherwise. Suppose $\beta \neq 0$. Then
%${\mathcal V}_{\alpha,\beta}= \{0\} \cup \{\gamma \in \Z_p : \gamma^{\ell} = (-\alpha \beta^{-1})^\frac{\ell}{d-1}\}$. 
%For each $\gamma \in \Z_p$ satisfying $\theta=\gamma^{\ell}$, for some $\theta \in \Z_p$, we have
%$\theta = (\gamma \omega^{\frac{p-1}{\ell}})^{\ell}$ where $\omega$ is a generator of $\Z_p^\times$.
%Therefore, for $\alpha \neq 0$ we have $|\{x \in \Z_p : x^{\ell} = (-\alpha \beta^{-1})^\frac{\ell}{d-1}\}| \in \{0, \ell \}$ and there are $\frac{p-1}{\ell}$ distinct $\ell$th roots of unity in $\Z_p^\times$. Since for each $\beta$, the map $\alpha \mapsto -(\alpha \beta^{-1})^\frac{\ell}{d-1}$ is a permutation, there are a total of $\frac{(p-1)^2}{\ell}$ pairs $\alpha, \beta \in \Z_p^\times$ such that $(-\alpha \beta^{-1})^\frac{\ell}{d-1}$ has an $\ell$th root in $\Z_p^\times$, in which case $c^f_{\alpha,\beta} = p-|{\mathcal V}_{\alpha,\beta}| =p-\ell-1$. 
 %\qed
 
 %$C_{x^d}$ defined over $\Z_{2p}$ is also a two-weight code.

We compute the weight distribution of a family of cyclic codes $C_{x^d}$ on $\Z_{2p}$,
which turn out to be two-weight codes. 

\begin{theorem}
    Let $R=\Z_{2p},p$ prime. Let $d \in \{2,...,p-1\}$, let $\ell=(d-1,p-1)$
    and let
   $f:R \longrightarrow R : x \mapsto x^d. $
    Then $$\Lambda_{f} = \{2p,\frac{2p\ell}{p-1},0 \}.$$ 
    and $C_{x^d}$ has weight enumerator
    $$1+ \frac{(p-1)^2}{\ell}X^{2p\frac{p-1-\ell}{p-1}}+ \left(2p^2- \frac{(p-1)^2}{\ell}-1\right)X^{2p}.$$
 \end{theorem}
 
 {\bf Proof:}
 The map $x \mapsto \alpha x+\beta x^d$ is identically zero if and only if $\alpha = \beta \in pR$. 
In particular, $|C_f| =2p^2$ and $c^f_{\alpha,\beta} = c^f_{\alpha',\beta'}$ if and only if $\alpha = \alpha' + p, \beta = \beta' + p$. 
We compute
\begin{eqnarray*}
W^{f}(\alpha,\beta)& = & \frac{1}{p-1}\sum_{x \in R} \left(\sum_{u \in R}\chi(u(\alpha x+ \beta x^d)) \right.
                             - \sum_{u \in pR} \chi(u(\alpha x+ \beta x^d))\\
                   & - & \left.\sum_{u \in 2R} \chi(u(\alpha x+ \beta x^d))  +1\right)\\
                   & = & \frac{1}{p-1}(2p|X| -2|X_2| -p |X_p| +2p),
\end{eqnarray*}
where $X = \{ x \in R : \alpha x + \beta x^d = 0\}$, $X_2 = \{ x \in R : \alpha x + \beta x^d \in 2R\}$ and $X_p = \{ x \in R : \alpha x + \beta x^d \in pR\}$.
%As in the proof of Theorem \ref{thzp}, 
For each $\gamma \in R^\times$, the map $\alpha \mapsto (-\alpha \gamma^{-1})^{\frac{\ell}{d-1}}$ is a permutation of $2R\backslash \{0\}$ if $\alpha \in 2R\backslash\{0\}$ and of $R^\times$ if $\alpha \in R^\times$. Moreover, for each divisor $\ell$ of $p-1$, there are $\frac{p-1}{\ell}$ distinct $\ell$th roots of unity in $R^\times$ and hence $\frac{(p-1)^2}{\ell}$ pairs $\alpha,\gamma$ such that $(-\alpha \gamma^{-1})^\frac{\ell}{d-1}$ has an $\ell$th root in $R^\times$ if $\alpha \in R^\times$ (resp. $2R $ if $\alpha \in \backslash \{ 0 \}$) . The results are summarized below.
\begin{center}
\begin{tabular}{|c|c|c|c|l|}
\hline
$|X|$ & $|X_2|$ & $|X_p|$ & $W^f(\alpha,\beta)$& $\alpha,\beta$ \\
\hline
$2(\ell+1)$ & $2p$ & $2(\ell+1)$ & $\displaystyle{\frac{2p \ell}{p-1}}$ & $\alpha=2\alpha_1, \beta=2\beta_1 \in 2R \backslash \{0\}$,\\ 
 & & & & $-\alpha_1\beta_1^{-1}$ has a $v$th root in $R$ \\
\hline
$2$ & $2p$ & $2$ & $0$ & $\alpha=2\alpha_1, \beta=2\beta_1 \in 2R \backslash \{0\},$\\ 
 & & & & $ -\alpha_1\beta_1^{-1}$ has no $v$th root in $R$\\
\hline
$\ell+1$ & $p $ & $2(\ell+1) $ & $0 $ & $\alpha = 2\alpha_1 \in 2R \backslash \{0\},$\\ 
 & & & & $\beta \in R^\times -2\alpha_1\beta^{-1}$ has a $v$th root in $2R$ \\ 
\hline
$1 $ & $p $ & $2 $ & $0 $ & $\alpha = 2\alpha_1 \in 2R \backslash \{0\},\beta \in R^{\times},$\\ 
 & & & & $ -2\alpha_1\beta^{-1}$ has no $v$th root in $2R$ \\ 
\hline
$2 $ & $2p $ & $2 $ & $0 $ & $\alpha=0,\beta \in 2R \backslash \{0\}$ \\ 
\hline
$1 $ & $p $ & $2 $ & $0 $ & $\alpha = 0, \beta \in R^{\times}$ \\ 
\hline
$ 2$ & $2p $ & $2 $ & $0 $ & $\alpha \in 2R\backslash \{0\},\beta = 0$ \\ 
\hline
$ 2$ & $2p $ & $2 $ & $0 $ & $\alpha \in R^{\times}, \beta = p $ \\ 
\hline
$p $ & $p $ & $2p $ & $0 $ & $\alpha =p, \beta = 0 $ \\ 
\hline
$2p $ & $p $ & $2 $ & $2p $ & $ \alpha = \beta = 0$ \\ 
\hline
\end{tabular}
\end{center}
~\qed

Given a two-weight code $C$ with non-zero weights $w_1<w_2$, we form the graph $G(C)=(V,E)$ by setting $V=C$ and $(x,y) \in E$ if and only if $w(x-y) = w_1$. A linear code $C$, with alphabet a finite Frobenius ring $R$ and generated by the $k \times n$ matrix $Y=(y_1,...,y_n)$ over $R$ is called {\em modular} if there is some $r \in \Q$ such that $|\{i: y_iR = y_jR\}| = r |y_jR^\times|$ for each $j$.
\begin{theorem}[Honold \cite{h07}]\label{thhon}
    Let $C$ be a modular linear code defined over a finite Frobenius ring with exactly two distinct non-zero weights. Then $G(C)$ is a strongly regular graph.
\end{theorem}
In fact for the finite field case, Theorem \ref{thhon} is easily deduced from \cite[Th. 2]{D72}. We now show that the code $C_{x^d}$ on $\Z_p$ determines a strongly regular graph.

\begin{corollary}
  Let $R=\Z_p$, $p$ prime, let $d \in \{2,...,p-1\}$ and let $\ell=(d-1,p-1)$. Then $C_{x^d}$ determines a strongly regular graph.% with parameters 
  %$$[p^2,\frac{(p-1)^2}{\ell},\frac{(p-1)^2}{\ell} - 1 -(p-1-\ell)(\ell(p+1)-p+1),(p-1)(p-1-\ell)].$$ 
\end{corollary} 

{\bf Proof:}
To establish the modular property we must show that $|\{y \in \Z_p^{\times} : c_{\alpha,\beta}^{x^d}(y) = \lambda c_{\alpha,\beta}^{x^d}(x)\; \forall \alpha,\beta \in \Z_p \}| = r(p-1)$, independently of our choice of $x$ for some $r \in \Q$. 
Let $x \in \Z_p^\times$. Then 
$$\lambda c_{\alpha,\beta}^{x^d}(x) = \alpha (\lambda x) + \beta \lambda^{p-d}(\lambda x)^d = 
 c_{\alpha,\beta \lambda^{p-d} }^{x^d}(\lambda x)=  c_{\alpha,\beta }^{x^d}(\lambda x),$$
 for all $\alpha,\beta \in \Z_p$ if and only if $\lambda^{d-1} = 1$.
 Then 
 $$|\{y \in \Z_p^{\times} : c_{\alpha,\beta}^{x^d}(y) = \lambda c_{\alpha,\beta}^{x^d}(x)\; \forall \alpha,\beta \in \Z_p \}|=| \{ \lambda x : \lambda^{d-1}=1\}| = \ell .$$
Thus $C_{x^d}$ is a modular two-weight code. The parameters $[n,k,\lambda,\mu]$ are uniquely determined by the parameters of $C_{x^d}$ (c.f. \cite{bgh,h07})
\qed

\begin{remark}
The codewords of $C_{x^d}$ defined over $\Z_p$ form an orthogonal array $OA(p,p-1)$ whenever $(d-1,p-1)=1$.
It is not hard to see that for $(d-1,p-1)=\ell >1$, the graph $G(C)$ is isomorphic to one induced by
an orthogonal array $OA(p,\frac{p-1}{\ell})$.
\end{remark}

\begin{remark}
The codes $C_{p\pi(x_0)x_1,\Z_{p^2}}$ and the codes $C_{x^d}$, defined on $\Z_{2p}$, are never modular, so that Theorem \ref{thhon} does not apply even though they are two-weight codes. Indeed in general the resulting graphs 
$G(C_{p\pi(x_0)x_1,\Z_{p^2}})$ and $G(C_{x^d})$ are disconnected.
\end{remark}

%
%{}
%Let $\Omega \subset R^k \backslash \{ 0 \}$ satisfy $\Omega R^\times = \Omega$. Let $M < {_RR^k}$.
%For each $v \in M$, define a function $c^v:\Omega \longrightarrow R: x \mapsto \langle v,x \rangle$.
%The left $R$-linear code $C:=C(\Omega)=\{c^v: v \in M \}$. 
%The Cayley graph of $\Omega$ is the graph $G:=G(\Omega)$ with vertex set $M$, where $u,v \in M$ are adjacent in
%$G(\Omega)$ if $u-v \in \Omega$. 
 %\begin{theorem}     
 %\end{theorem}
%

\subsection{Local Commutative Rings}

For this construction, we assume the existence of $\sigma \in Aut(R)$ such that
$\chi \circ \sigma = \chi $ on $R$. For $R$ a finite field, the function described in the next lemma is
the map $:x \mapsto x^{p^k+1}$, which is known to be perfect nonlinear on a field of odd characteristic, and almost perfect nonlinear on a field of even characteristic.
 
\begin{lemma}\label{lemfaut}
   Let $\sigma \in Aut(R)$ satisfy $\chi(\sigma (x)) = \chi(x)$ for 
   all $x \in R$.
   Define 
   $$f: R \longrightarrow R : a \mapsto \sigma(a)a - \sigma(a_m)a_m. $$
   Then 
   $$W^{f,S}(\alpha,\beta) = \frac{|M|}{|S^\times|}\left( |S||{\mathcal V} \cap U_0| - |S \cap M||{\mathcal V} \cap U_1|   \right),$$
   where ${\mathcal V} = \{ t \in {\mathcal T} : \alpha + \beta \sigma(t) + \sigma^{-1}(\beta t) \in \soc R\}$,
   $U_0= \{ t \in {\mathcal T} : T(\alpha t + \beta \sigma(t) t)  =0 \}$ and $U_1= \{ t \in {\mathcal T} : T(\alpha t + \beta \sigma(t) t ) \in \soc S\}$.
   \end{lemma} 
   
   {\bf Proof:} Let $\alpha, \beta \in R$
   \begin{eqnarray*}
      W^{f,S}(\alpha,\beta) & = & \frac{1}{|S^\times|} \sum_{u \in S^\times} \sum_{x \in R} \chi^u(\alpha x + \beta (\sigma(x)x -\sigma(x_m)x_m) ),\\
      & = & \frac{1}{|S^\times|} \sum_{u \in S^\times} \sum_{t \in {\mathcal T}}  \sum_{m \in M} \chi^u(\alpha(t + m) + \beta (\sigma(t+m)(t+m) -\sigma(m)m) ), \\
      & = & \frac{1}{|S^\times|} \sum_{u \in S^\times} \sum_{t \in {\mathcal T}} \chi^u(\alpha t + \beta \sigma(t)t)
       \sum_{m \in M} \chi^u(\alpha m + \beta(\sigma(t)m + \sigma(m) t)),\\
       & = & \frac{1}{|S^\times|} \sum_{u \in S^\times} \sum_{t \in {\mathcal T}} \chi^u(\alpha t + \beta \sigma(t)t)
       \sum_{m \in M} \chi^u((\alpha + \beta \sigma(t)+ \sigma^{-1}(\beta t))m),
   \end{eqnarray*}
   since $\chi(\sigma(x))=\chi(x)$ for all $x \in R$. The character $\chi^u((\alpha + \beta \sigma(t)+ \sigma^{-1}(\beta t))\cdot)$ is trivial on $M$ if and only if $\alpha + \beta \sigma(t)+ \sigma^{-1}(\beta t)$ annihilates every element of $M$. Therefore
   $$W^{f,S}(\alpha,\beta) =  \frac{|M|}{|S^\times|} \sum_{u \in S^\times} \sum_{t \in {\mathcal V}} \chi^u(\alpha t + \beta \sigma(t)t),$$
   where ${\mathcal V} = \{ t \in {\mathcal T} : \alpha + \beta \sigma(t) + \sigma^{-1}(\beta t) \in \soc R\}$. 
   Since $S$ is local, we have $S^\times = S \backslash S \cap M$ and the annihilator ideal of $S\cap M$ in $S$ is $\soc S$, from which we deduce
   \begin{eqnarray*}
   W^{f,S}(\alpha,\beta)  & = &  \frac{|M|}{|S^\times|} \sum_{t \in {\mathcal V}}
   \left( \sum_{u \in S} \Phi(uT(\alpha t + \beta \sigma(t)t)) - \sum_{u \in S \cap M} \Phi(uT(\alpha t + \beta \sigma(t)t) )\right)\\
   & = & \frac{|M|}{|S^\times|} \left( |S||{\mathcal V} \cap U_0| - |S \cap M||{\mathcal V} \cap U_1|   \right),
   \end{eqnarray*}
   where
   $U_0= \{ t \in {\mathcal T} : T(\alpha t + \beta \sigma(t) t)  =0 \}$ and $U_1= \{ t \in {\mathcal T} : T(\alpha t + \beta \sigma(t) t ) \in \soc S\}$.
   
\begin{example}\label{exlocfrob}
   Let $R = \Z_2[x,y] / \langle x^2,y^2 \rangle = \{a_1+a_x x+a_y y +a_{xy} xy: a_X \in \Z_2 \}$. 
   $R$ is a finite local Frobenius ring with Teichmuller set $\Z_2$ and maximal ideal 
   $\langle x, y \rangle$. Each $a \in R$ has a unique expression in the form $a = a_t+a_m$ with $a_t = a_1 
   \in \Z_2$ and $a_m = a_y x + a_x y + a_{xy} xy$.  
   The automorphism group of $R$ has order 2 and is generated by 
   $$\sigma :R \longrightarrow R : a_1+a_x x + a_y y + a_{xy} xy \mapsto a_1+a_y x + a_x y + a_{xy} xy.$$
   It can be checked that the map 
   $$T : R \longrightarrow \Z_2 : a_1+a_x x + a_y y + a_{xy} xy \mapsto a_1 +a_x+a_y+a_{xy}$$
   is a trace map from $R$ onto $\Z_2$ that induces the character
   $\chi:R \longrightarrow \C^\times$ defined by $\chi(a) = (-1)^{T(a)}.$ Moreover, as $T\circ \sigma=T$,
   $\chi(\sigma(a))=\chi(a)$ for all $a\in R$.
   If $f:R \longrightarrow R$ is defined as in Lemma \ref{lemfaut} then
   $f(a) =a_1(1+(a_x+a_y)(x+y))$ for each $a\in R$.
   With the same notation as in the lemma, for any $\alpha,\beta \in \R$, we have $U_0= \{t \in\Z_2: tT(\alpha+\beta) = 0 \}$ and $U_1 = \Z_2$ 
   so that $W^{f,\Z_2}(\alpha,\beta) = 8(2|{\mathcal V} \cap U_0|-|{\mathcal V}|)$, where 
   ${\mathcal V} = \{t \in \Z_2: \alpha + (\beta + \sigma(\beta))t \in \langle xy \rangle \}$.
   %Note that %${\mathcal V}$ is non-empty if and only if $\alpha \in \langle x+y \rangle$ and $a_x=a_y=t(\beta_x+\beta_y)$.
   %$$ {\mathcal V} =\left\{ \begin{array}{cl}
    %                        2 & \mbox{ if } \alpha \in \langle xy \rangle \mbox{ and } \beta_x = \beta_y\\
    %                       %0 & \mbox{ if } \alpha \notin \langle x+y \rangle\\
    %                        1 & \mbox{ if } \alpha \in  \langle x+y \rangle \mbox{ and } \beta_x\neq\beta_y\\
    %                        0 & \mbox{otherwise}
    %                        \end{array} \right. $$
   It is straightforward to check that $\Lambda_{f,\Z_2} = \{-8,0,8,16\}$ and that $C_{f,\Z_2}$ is a 
   three-weight code with Hamming weights $\{0,4,8,12\}$ (setting $\gamma = \frac{q-1}{q}$ in the homogeneous weight).
   Its Hamming weight enumerator is given by $1+3X^4+27X^{8}+X^{12}$.
\end{example}

\begin{example}\label{exgallocfrob}
  Let $R$ be the Galois ring GR$(2^3,2)$ of characteristic $8$ and order $64$.
  As we saw in Example \ref{exgaltr}, $Aut(R)$ is cyclic of order $2$, generated by
  $\sigma$ satisfying $T_{R/\Z_8}= T_{R/\Z_8}\circ \sigma$. Then with $f$ defined as in 
  Lemma \ref{lemfaut}, we have 
  $$f(a_0+2a_1+4a_2) = 1 + 2(1+a_0^2a_1+a_1^2a_0)+4(1+a_0^2a_2+a_2^2a_0).$$
  We compute the distinct homogeneous weights of code $C_{f,\Z_8}$ as 
  $\{0,32,48,64,80,88,96 \}$.
\end{example}

\begin{corollary}\label{corlocfrob}
   Let  $f: R \longrightarrow R$ be defined as in Lemma \ref{lemfaut}
   Then $$\Lambda_{f,R} =\{ |R|,|M|,\frac{|R||M|}{|R^\times|} ,0  \}.$$
   If $|K|>2$ then $C_{f,R} $ has size $|R|^2$  and weight distribution    
   \begin{center}
   \begin{tabular}{|c|c|}
   \hline
   weight & number of codewords\\
   \hline
       $0$ & $1$ \\
       \hline
       $|R|-|M|$ & $|K||R| - |K|^2 + |K| - 1$ \\
       \hline
       $\displaystyle{|R| - \frac{|R||M|}{|R^\times|}}$ & $|K|^2-2|K| +1$ \\
       \hline
       $|R|$ & $|R|^2 - |K||R| + |K| - 1$\\
   \hline
   \end{tabular}.
   \end{center}
   If $|K|=2$ then $\displaystyle{|C_{f,R}|=\frac{|R|^2}{2}}$ and $C_{f,R}$ has weight enumerator
   $$1+ (|R|-2)X^{\frac{|R|}{2}}+\left(\frac{|R|^2}{2}-|R|+1\right)X^{|R|}. $$
\end{corollary} 

{\bf Proof:}
   We note that every element of $Aut(R)$ fixes $R^\times$, $M$ and $\soc R$.
   Let $\alpha,\beta \in R$. 
   If $\alpha,\beta \notin \soc R$ then $0 \notin {\mathcal V}$ and so  
   $${\mathcal V} \cap U_1 = \{t \in {\mathcal T}\backslash \{0\} : \sigma^{-1}(\beta t) \in \soc R \} = \emptyset = {\mathcal V} \cap U_0 .$$ 
   
   If $\alpha \in \soc R, \beta \notin \soc R$ then  
   $$\{0\} \subset {\mathcal V} \cap U_1 \subset \{t \in {\mathcal T} : \beta \sigma(t)t \in 
   \soc R \} = \{ 0\} = {\mathcal V} \cap U_0.$$ 
   
   If $\alpha \notin \soc R, \beta \in \soc R$ then ${\mathcal V} = \emptyset$. 
   
   Suppose now that $\alpha,\beta \in \soc R$. Then clearly ${\mathcal V} = U_1 = {\mathcal T}$. 
   Since $\soc R=M^\perp$ is principal in $R$, there exist unique $\alpha_u,\beta_u \in {\mathcal T}$ satisfying 
   $\alpha=\alpha_u \theta$ and $\beta=\beta_u \theta$ for fixed $\theta$ generating $\soc R$. Therefore,
   $$U_0=\{0\} \cup \{ t \in {\mathcal T} \backslash \{0\}: \alpha_u + \beta_u \sigma(t) \in M \}.$$ 
   If $\beta_u \neq 0$ then there is a unique $t \in {\mathcal 
   T}$ satisfying $t \in \sigma^{-1}(-\alpha_u \beta_u^{-1})+M$. If $\alpha_u =0 $ we have $|U_0| = 1$ and otherwise we have $|U_0| = 2$. If 
   $\beta_u=0$ then $U_0 = \{ 0 \}$ unless $\alpha_u =0$, in which case $U_0 = {\mathcal T}$.
   We summarize these observations in the table below.
   
   \begin{center}
   \begin{tabular}{|c|l|}
   \hline
   $W^f(\alpha,\beta)$ & $\alpha,\beta$\\
   \hline
       $|R|$ & $\alpha=\beta=0$ \\
       \hline
       $|M|$ & $\alpha \in \soc R, \beta \notin \soc R$ or $\alpha \in \soc R \backslash \{0\},\beta=0$ \\
       \hline
       $\displaystyle{\frac{|R||M|}{|R^\times|}}$ & $\alpha,\beta \in \soc R \backslash \{0\}$ \\
       \hline
       $0$ & $\alpha \notin \soc R$ or $\alpha=0, \beta \in \soc R \backslash \{0\}$\\
   \hline
   \end{tabular}.
   \end{center}
   
   Now $c^f_{\alpha,\beta}(x) =\alpha x + \beta(\sigma(x_t)x+\sigma(x_m)x_t)$.
   Suppose that $c^f_{\alpha,\beta}$ is identically zero. Then in particular, 
   $c^f_{\alpha,\beta}(\theta) = \alpha \theta + \beta(\sigma(\theta)\theta = 0$ for any $\theta \in{\mathcal T},$
   in which case $ \alpha  = - \beta \sigma(\theta)$ for any such nonzero $\theta$. It follows that if $|K|=|{\mathcal T}|>2$
   then $\alpha=\beta=0$. Therefore $|C_{f,R}|=|R|^2$, and $C_{f,R}$ has the weight distribution shown above.
   Suppose now that $|K|=|{\mathcal T}|=2$, so that ${\mathcal T}=\{0,1\}$, which is point-wise fixed by $\sigma$. Then 
   $\alpha = -\beta$ and moreover, since $c^f_{\alpha,\beta}(m) = \alpha m$ for any $m\in M$, we must have $\alpha \in \soc R$.
   Conversely, for any $\delta \in \soc R$,  
   $c^f_{\delta,-\delta}(x)=\delta(x(1-x_t)+\sigma(x_m)x_t)$ is the zero map.
   It follows that $\displaystyle{|C|=\frac{|R|^2}{|\soc R|} = \frac{|R|^2}{2}}$ and that $C_{f,R}$ has the weight enumerator given above.  
   
\qed

\begin{example}
   Let $R = \Z_2[x,y] / \langle x^2,y^2 \rangle$ with $f$ defined as in Example \ref{exlocfrob}. Since $|K|=2$,
   from Corollary \ref{corlocfrob} we deduce that $C_f$ is a two-weight code over $R$ of size 128 with homogeneous weight enumerator
   $$1+ 14 X^{8} +113 X^{16}.$$  
\end{example}

\begin{example}
   Let $R = \Z_3[x,y] / \langle x^2,y^2 \rangle$. Then $R$ has residue field $K=GF(3)$, $|R|=3^4,|M|=3^3$ and $|R^\times|=2.3^3$. 
   Let $f(x) = \sigma(x)x-\sigma(x_m)x_m$ for 
   $\sigma(a_1+a_xx+a_yy+a_{xy}xy)=a_1+a_yx+a_xy+a_{xy}xy$. Then from Corollary \ref{corlocfrob}, $C_f$ is a three-weight code over $R$ of size
   $3^8=6561$ with weight enumerator
   $$1+4 X^{\frac{81}{2}}+236 X^{54} + 6320 X^{81}.$$   
\end{example}

\begin{example}
   Let $R=\mathrm{GR}(2^3,2)$, with $f$ defined as in Example \ref{exgallocfrob}.
   Then $C_f$ is a three-weight code over $R$ of order $4096$ with weight enumerator
   $$1 + 243 X^{48}+ 9X^{\frac{128}{3}} + 3843 X^{64}.$$
\end{example}

%
%{Local Commutative Rings}

%\begin{theorem}
 %  Let $R$ be a finite local commutative Frobenius ring. Let $S$ be a subring of $R$ such that there exists
 %  a trace map $T:R \longrightarrow S$. 
  % Let $\sigma \in Aut(R)$ satisfy $\chi(\sigma (x)) = \chi(x)$ for %all $x \in R$.
 %  Define $$f: R \longrightarrow R : a \mapsto \sigma(a)a - \sigma(a_m)a_m. $$
  % Then $$\Lambda_f =\{ |R|, \frac{|M|}{|S^\times|}(|S| - |M_S||T|) ,0  \}. $$
%\end{theorem} 
%

%
%{Local Commutative Rings}

%\begin{theorem}
 %  Let $R$ be a finite local commutative Frobenius ring. Let $S$ be a subring of $R$ for which there exists
  % a trace map $T:R \longrightarrow S$. Let $\sigma \in Aut(R)$ satisfy $\chi(\sigma (x)) = \chi(x)$ for all $x \in R$.
   %Define $$f: R \longrightarrow R : a \mapsto \sigma(a)a - \sigma(m_a)m_a. $$
   %Then $C_f $ has non-zero weights
   %$$w_1= |R|, w_2=|R|+\frac{|M|}{|S^\times|}( |M_S||T|-|S|). $$
%\end{theorem} 
%

{\bf Acknowledgement: } The author is grateful to Marcus Greferath and Alexander Nechaev for useful discussions regarding trace maps on finite rings.


\begin{thebibliography}{99}

\bibitem{bbmmg} 
 C. Bracken, E. Byrne, N. Markin, G. McGuire \emph{New Families of Almost Perfect Nonlinear Trinomials and Multinomials}, Finite Fields and Their Applications, 14 (3) (2008) 703--714.

\bibitem{bgh}
E. Byrne, M. Greferath and T. Honold, \emph{Ring Geometries, Two-Weight Codes and Strongly Regular Graphs}, Designs, Codes and Cryptography, 48 (1) (2008) 1--16.

\bibitem{bgos}
E. Byrne, M.~Greferath and M.~E.~O'Sullivan, {\em The Linear Programming Bound for Codes 
over Finite Frobenius Rings,} Designs, Codes and Cryptography,
42 ,  {\textbf 3}  (2007), 289--301.

\bibitem{bgks}
E. Byrne, M. Greferath, A. Kohnert, V. Skachek, {\em New Bounds for Codes Over Finite Frobenius Rings}
Designs, Codes and Cryptography, 42 {\textbf 3}  (2010), Online First.

\bibitem{bs}
E. Byrne, A. Sneyd, \emph{Constructions of Two-Weight Codes Over Finite Rings},
Proceedings of the 19th International Symposium on Mathematical Theory of Networks and Systems (MTNS 2010),
Budapest, July, 2010.

\bibitem{CCZ} C. Carlet, P. Charpin, V. Zinoviev, {\em Codes, Bent Functions and Permutations
Suitable for DES-Like Cryptosystems}, Designs, Codes and Cryptography,
Vol. 15, No. 2, (1998) 125--156.

\bibitem{cdy}
C. Carlet, C. Ding, J. Yuan, \emph{Linear Codes From Perfect Nonlinear Mappings and Their Secret Sharing Schemes}, IEEE Trans. Inform. Th., Vol. 51, {\bf 6}, (2005) 2089--2013

\bibitem{ch}
I.~Constantinescu and W.~Heise, \emph{A Metric for Codes over Residue Class
Rings of Integers}, Problemy Peredachi Informatsii, 33 (3) (1997).

\bibitem{D72} P. Delsarte, \emph{Weights of linear codes and strongly regular normed spaces}, Discrete Math., {\bf 3} (1972) 47--64. 
  
\bibitem{bruce}
M.~Greferath and S.~E.~Schmidt, \emph{Finite-Ring Combinatorics and MacWilliams Equivalence Theorem}, J. of Combinatorial Theory (A) \textbf{92} (2000), 17--28. 

\bibitem{gnw}
M.~Greferath, A. Nechaev, R. Wisbauer, \emph{Finite Quasi-Frobenius Modules and Linear Codes}, Journal of Algebra and its Application, {\bf 3} (3) (2004) 247--272.

\bibitem{plotkin}
M.~Greferath and M.~E.~O'Sullivan, \emph{On Bounds for Codes over Frobenius Rings under Homogeneous Weights}, Discrete Mathematics \textbf{289} (2004), 11--24.

  
\bibitem{hammons94}
{A. R.} Hammons, {P. V.} Kumar, {A. R.} Calderbank, {N. J. A. } Sloane, and
{P.} Sol\'e, \emph{The ${\mathbb Z}_4$-linearity of {K}erdock, {P}reparata,
{G}oethals and related codes}, IEEE Trans. Inform. Theory \textbf{40} (1994),
301--319.  
  
\bibitem{frank}
R. C. Heimiller, \emph{Phase Shift Pulse Codes with Good Periodic Correlation Properties}, IRE Transactions on Information Theory, IT-7 (1961) 254--257.

\bibitem{h}
T. Honold, \emph{A Characterization of Finite Frobenius Rings}, Arch. Math. (Basel), {\bf 76} (2001).
 
\bibitem{h07}
T. Honold, {\em Further Results on Homogeneous Two-Weight Codes}, Proceedings of
Optimal Codes and Related Topics, Bulgaria (2007).

\bibitem{lam} T. Y. Lam, \emph{Lectures on Modules and Rings},
Graduate Texts in Mathematics, Vol. 189, Springer-Verlag, 1999.

\bibitem{macd}
B. R.~McDonald, \emph{Finite Rings With Identity}, Pure and Applied Mathematics, M. Dekker, New York (1974).

\bibitem{r}
R. Raghavendran, {\em Finite Associative Rings}, Compositio Math., 21 (1969) 195-229.
%\bibitem{wood97a}
%{J. A.} Wood, \emph{Duality for modules over finite rings and applications to
%coding theory}, Amer. J. Math. \textbf{121} (1999), 555--575. 

\bibitem{ycd}
J. Yuan, C. Carlet, C. Ding, \emph{The Weight Distribution of a Class of Linear Codes From Perfect Nonlinear Functions}, IEEE Trans. Inform. Th., Vol. 52, {\bf 2} (2006) 712--717.

 
 \end{thebibliography}
 \end{document}